\documentclass{amsart}

\usepackage{amsmath,amssymb,amscd}
\usepackage{epsfig,pstricks}
\usepackage{xypic}

\newtheorem{theorem}{Theorem}[section]

\newtheorem{lemma}[theorem]{Lemma}

\theoremstyle{definition}
\newtheorem{remark}[theorem]{Remark}
\newtheorem{example}[theorem]{Example}
\newtheorem{definition}[theorem]{Definition}

\renewcommand{\leq}{\leqslant}
\renewcommand{\geq}{\geqslant}

\newcommand{\assign}{:=}
\newcommand{\colim}{{\rm colim}}
\renewcommand{\lim}{{\rm lim}}

\newcommand{\xra}[1]{\xrightarrow{#1}}

\newcommand{\lmod}[1]{{\rm {\bf Mod}}({#1})}
\newcommand{\comod}[1]{{\rm {\bf CoMod}}({#1})}

\newcommand{\C}[1]{\mathcal{#1}}

\newcommand{\B}[1]{\mathbb{#1}}
\newcommand{\T}[1]{{\tt #1}}

\newcommand{\app}[2]{{\rm App}_{#1}(#2)}

\numberwithin{equation}{section}

\title{The universal Hopf cyclic theory}
\author{Atabey Kaygun}
\address{Department of mathematics, The University of Western Ontario, London N6A 5B7 Canada}
\email{\tt akaygun@uwo.ca}

\begin{document}
\maketitle

\section{Introduction}

In noncommutative geometry, as the category of algebras of various
flavors replaced the category of spaces of various flavors, Hopf
algebras arose as the natural candidate to study the symmetries of a
noncommutative space.  Unlike the classical notion of symmetry, the
notion of noncommutative symmetry has four different types:
\begin{quote}
\begin{tabular}{lll}
 {\tt [MC]} &   {\bf module coalgebra}   & A Hopf algebra acting on a coalgebra in a compatible way.\\ 
 {\tt [CA]} &   {\bf comodule algebra}   & A Hopf algebra coacting on a algebra in a compatible way.\\ 
 {\tt [MA]} &   {\bf module algebra}     & A Hopf algebra acting on a algebra in a compatible way.\\
 {\tt [CC]} &   {\bf comodule coalgebra} & A Hopf algebra coacting on a coalgebra in a compatible way.
\end{tabular}
\end{quote}
These compatibility conditions can be expressed concisely as the
(co)multiplication structure morphism of the corresponding (co)algebra
being a $B$-(co)module morphism where $B$ is our base Hopf algebra.
We are interested in such symmetries in the context of cyclic
(co)homology.  In the sequel the term ``cyclic theory'' will mean a
functor from a suitable category of algebras into the category of
(co)cyclic $k$-modules and the term ``cyclic (co)homology'' will mean
a suitable (co)homology functor from the category of (co)cyclic
$k$-modules into the category of $k$-modules.

Since Connes and Moscovici's seminal work on the cyclic cohomology of
Hopf algebras and transverse index
theorem~\cite{ConnesMoscovici:HopfCyclicCohomology,
ConnesMoscovici:HopfCyclicCohomologyII}, there has been a surge in the
interest in such symmetries
\cite{Crainic:CyclicCohomologyOfHopfAlgebras, JaraStefan:HopfGalois,
Menichi:BatalinVilkoviskyAlgebras,
NeshveyevTuset:HopfEquivariantCyclicCohomology,
Sharygin:CyclicCohomologyViaKaroubi,
Taillefer:CyclicHomologyOfHopfAlgebras}.  However, a (co)cyclic theory
and basic tools of cyclic cohomology had to be built from scratch for
each type of symmetry separately~\cite{Khalkhali:HopfCyclicHomology,
Kaygun:BialgebraCyclicK, Kaygun:BivariantHopf}.  To complicate the
story further, the cyclic dual of each (co)cyclic theory constructed
thus far is non-trivial in stark contrast with the ordinary
non-equivariant case where the dual is always trivial.  This means
there are 8 potentially different types of cyclic theories in the
presence of a Hopf symmetry.  Considering the variety of (co)homology
functors one can apply to a given (co)cyclic object, we suddenly see a
phletora of cyclic cohomology theories which are sensitive to the Hopf
symmetry in the noncommutative universe.
 
Amazingly, Khalkhali and Rangipour~\cite{Khalkhali:CyclicDuality}
showed that if we view a Hopf algebra as a coalgebra enjoying an
\T{[MC]}-type symmetry over itself, then the cyclic dual of the
canonical cocyclic object is functorially isomorphic to the canonical
cyclic object of the same Hopf algebra viewed as an algebra enjoying
an \T{[CA]}-type symmetry over itself.  This result suggests that
there is a deep meta-symmetry lurking behind, connecting all of these
different types of cyclic theories.

In a previous paper~\cite{Kaygun:BivariantHopf}, Khalkhali and the
author successfully unified the \T{[MA]} and \T{[MC]}-type cyclic
theories and their cyclic duals under the banner of bivariant
Hopf-cyclic cohomology.  In this paper, we aim to unravel further the
meta-symmetry behind all of these cyclic theories and construct a new
universal cyclic theory encompassing all types of symmetries we stated
above, agreeing with, and extending further, the previous definitions
given in the literature.  The construction is purely categorical and
each individual theory is obtained by modifying certain parameters.
These parameters, or variables, are (i) a symmetric monoidal category
which will stand for the category modules over a ground ring $k$, (ii)
a class of morphisms called {\em transpositions} which will play the
role of a coefficient, (iii) an arbitrary exact comonad which will
replace a $k$-flat Hopf algebra and finally (iv) a suitable category
of (co)algebras called {\em transpositive (co)algebras} which will
play the role of (co)module (co)algebras.

One {\em practical} consequence of this formal exercise in category
theory is that we no longer need to define a different theory for each
type of symmetry and then prove directly that it really is cyclic,
which is quite technical and involved.  The recipe we provide in this
paper ensures that the end object is not only equivariantly (co)cyclic
but also the right object for all known cases.  The results of this
article will give us the licence to ignore the technical problems of
existence of a right kind of cyclic theory and to engage with more
pressing questions such as excision, Morita invariance and homotopy
invariance in the presence of a Hopf symmetry.  Moreover, now that the
these cyclic theories are defined by universal properties we expect
such questions to become more accessible for investigation.

Here is a plan of this paper.  In Section~\ref{sec:Transpositive} we
give definitions of transpositions and transpositive (co)algebras in
an arbitrary symmetric strict monoidal category $\C{C}$.  In the same
section we also describe ordinary $B$-(co)module (co)algebras over an
arbitrary bialgebra $B$ as transpositive algebras in a specific
monoidal category with respect to certain classes of transpositions.
In Section~\ref{sec:UniversalParaCocyclicTheory} we construct the
universal para-(co)cyclic theory for the category of transpositive
(co)algebras.  Next in Section~\ref{sec:Approximations}, we
incorporate an arbitrary exact comonad $\T{B}$ into our machinery.  In
this section, specifically in Theorem~\ref{thm:MainApproximation}, we
show that every pseudo-para-(co)cyclic $\T{B}$-comodule admits an {\em
approximation} (Definition~\ref{defn:Approximation}) in the category
of (co)cyclic $\T{B}$-comodules.  For an arbitrary bialgebra $B$, in
Section~\ref{sec:Examples} we recover the Hopf cyclic and equivariant
cyclic theories of $B$-module (co)algebras~\cite{Kaygun:BivariantHopf}
and bialgebra cyclic theory of $B$-comodule
algebras~\cite{Kaygun:BialgebraCyclicK}.  The key observation we use
is that the universal para-cyclic theory actually takes values in the
category of pseudo-para-(co)cyclic $B$-modules in these cases.  As a
side result, we recover the Hopf--Hochschild
homology~\cite{Kaygun:HopfHochschild} by using the techniques
developed in this paper.  We end the paper by defining the missing
cyclic theory for comodule coalgebras as a natural extension of the
cyclic theories defined hitherto.

Throughout this article, we assume $\C{C}$ is a small category.  If we
require $\C{C}$ to be monoidal $\otimes$ will denote the monoidal
product of $\C{C}$ and we will assume $(\C{C},\otimes)$ is a symmetric
strict monoidal category with a unit object $I$.

\subsection{Acknowledgements}

We trace the genesis of this article back to a series of discussions
the author had with Krzysztof Worytkiewicz on the formal properties of
the Hopf cyclic (co)homology.

\section{Transpositions and transpositive (co)algebras}
\label{sec:Transpositive}

In this section we will use a rudimentary version of ``calculus of
braid diagrams'' in monoidal categories as developed in
\cite{Majid:BraidCalculus}.

\begin{definition}
  Fix an object $M$ in $\C{C}$ and a unique morphism $w_{M,X}:M\otimes
  X\xra{}X\otimes M$ for each object $X$ chosen from a subset $\C{T}$
  of $Ob(\C{C})$.  The datum $(M,\C{T},\{w_{M,X}\}_{X\in\C{T}})$ is
  called {\em a class of transpositions} and is denoted by $w$.  For
  such an object $X\in\C{T}$, the morphism $w_{M,X}$ and its inverse
  $w_{M,X}^{-1}$ if it exists, are going to be denoted by
\begin{equation*}
\uput{0}[0](-2.60,-0.00){M}
\uput{0}[0](-2.10,-0.00){X}
\uput{0}[0](-0.10,-0.00){X}
\uput{0}[0](0.40,-0.00){M}
\psline(-2.5000,-0.2500)(-2.5000,-0.3750)
\psline(-2.0000,-0.3750)(-2.5000,-0.6250)
\psline[border=2pt](-2.5000,-0.3750)(-2.0000,-0.6250)
\psline(-2.0000,-0.6250)(-2.0000,-0.7500)
\psline(-2.0000,-0.2500)(-2.0000,-0.3750)
\psline(-2.5000,-0.6250)(-2.5000,-0.7500)
\psline(0.0000,-0.2500)(0.0000,-0.3750)
\psline(0.0000,-0.3750)(0.5000,-0.6250)
\psline[border=2pt](0.5000,-0.3750)(0.0000,-0.6250)
\psline(0.5000,-0.6250)(0.5000,-0.7500)
\psline(0.5000,-0.2500)(0.5000,-0.3750)
\psline(0.0000,-0.6250)(0.0000,-0.7500)
\uput{0}[0](-2.60,-1.00){X}
\uput{0}[0](-2.10,-1.00){M}
\uput{0}[0](-0.10,-1.00){M}
\uput{0}[0](0.40,-1.00){X}
\vspace{0.9500cm}
\end{equation*}
  respectively.  We do not require transpositions to be invertible,
  nor do we require the inverses to be transpositions themselves even
  if they exist.
\end{definition}

An object $A$ in $\C{C}$ is called an algebra if there exist morphisms
$A^{\otimes 2}\xra{\mu_A}A$ and $I\xra{e}A$ such that the following
diagrams commute
\begin{equation*}
  \xymatrix{
    \ar[d]_{A\otimes \mu_A} A^{\otimes 3} \ar[r]^{\mu_A\otimes A} & 
    A^{\otimes 2} \ar[d]^{\mu_A} 
            & & \ar[d]_{e\otimes I} A \ar[r]^{I\otimes e} \ar[rd]_{id} &  \ar[d]^{\mu_A} A^{\otimes 2}\\
    A^{\otimes 2} \ar[r]_{\mu_A}                 & 
    A       & &  A^{\otimes 2}  \ar[r]_{\mu_A}   &  A
  }
\end{equation*}
In other words, $A$ is a monoid object in $(\C{C},\otimes)$, or
equivalently both $(A\otimes\ \cdot\ )$ and $(\ \cdot\ \otimes A)$ are
a monads in $\C{C}$.  A coalgebra $(C,\delta_C,\varepsilon)$ in
$\C{C}$ is simply an algebra in $\C{C}^{op}$.  These conditions for an
algebra $A$ will also be denoted by the following diagrams
\begin{equation*}
\uput{0}[0](-4.35,-0.00){A}
\uput{0}[0](-3.85,-0.00){A}
\uput{0}[0](-3.35,-0.00){A}
\uput{0}[0](-2.35,-0.00){A}
\uput{0}[0](-1.85,-0.00){A}
\uput{0}[0](-1.35,-0.00){A}
\uput{0}[0](0.40,-0.00){A}
\uput{0}[0](1.40,-0.00){A}
\psline(-4.2500,-0.2500)(-4.2500,-0.3750)
\psline(-4.2500,-0.3750)(-4.0000,-0.6250)
\psline(-4.0000,-0.6250)(-4.0000,-0.7500)
\psline(-3.7500,-0.2500)(-3.7500,-0.3750)
\psline(-3.7500,-0.3750)(-4.0000,-0.6250)
\psline(-3.2500,-0.2500)(-3.2500,-0.7500)
\psline(-2.2500,-0.2500)(-2.2500,-0.7500)
\psline(-1.7500,-0.2500)(-1.7500,-0.3750)
\psline(-1.7500,-0.3750)(-1.5000,-0.6250)
\psline(-1.5000,-0.6250)(-1.5000,-0.7500)
\psline(-1.2500,-0.2500)(-1.2500,-0.3750)
\psline(-1.2500,-0.3750)(-1.5000,-0.6250)
\psline{o}(0.0000,-0.5000)(0.0000,-0.7500)
\psline(0.5000,-0.2500)(0.5000,-0.7500)
\psline(1.5000,-0.2500)(1.5000,-0.7500)
\psline{o}(2.0000,-0.5000)(2.0000,-0.7500)
\uput{0}[0](2.90,-0.50){A}
\psline(-4.0000,-0.7500)(-4.0000,-0.8750)
\psline(-4.0000,-0.8750)(-3.6250,-1.1250)
\psline(-3.6250,-1.1250)(-3.6250,-1.2500)
\psline(-3.2500,-0.7500)(-3.2500,-0.8750)
\psline(-3.2500,-0.8750)(-3.6250,-1.1250)
\uput{0}[0](-2.98,-1.00){~=}
\psline(-2.2500,-0.7500)(-2.2500,-0.8750)
\psline(-2.2500,-0.8750)(-1.8750,-1.1250)
\psline(-1.8750,-1.1250)(-1.8750,-1.2500)
\psline(-1.5000,-0.7500)(-1.5000,-0.8750)
\psline(-1.5000,-0.8750)(-1.8750,-1.1250)
\psline(0.0000,-0.7500)(0.0000,-0.8750)
\psline(0.0000,-0.8750)(0.2500,-1.1250)
\psline(0.2500,-1.1250)(0.2500,-1.2500)
\psline(0.5000,-0.7500)(0.5000,-0.8750)
\psline(0.5000,-0.8750)(0.2500,-1.1250)
\uput{0}[0](0.65,-1.00){~=}
\psline(1.5000,-0.7500)(1.5000,-0.8750)
\psline(1.5000,-0.8750)(1.7500,-1.1250)
\psline(1.7500,-1.1250)(1.7500,-1.2500)
\psline(2.0000,-0.7500)(2.0000,-0.8750)
\psline(2.0000,-0.8750)(1.7500,-1.1250)
\uput{0}[0](2.40,-1.00){=}
\psline(3.0000,-0.7500)(3.0000,-1.2500)
\uput{0}[0](-3.73,-1.50){A}
\uput{0}[0](-1.98,-1.50){A}
\uput{0}[0](0.15,-1.50){A}
\uput{0}[0](1.65,-1.50){A}
\uput{0}[0](2.90,-1.50){A}
\vspace{1.4500cm}
\end{equation*}

\begin{definition}
  An algebra $(A,\mu_A,e)$ in $\C{C}$ is called a {\em
  $w$-transpositive algebra} if there exists a morphism $w_{M,X}\colon
  M\otimes X\xra{}X\otimes M$ in $w$ such that the following diagram
  commutes
  \begin{equation*}
    \xymatrix{
      \ar[d]_{M\otimes\mu_A} M\otimes A\otimes A \ar[r]^{w_{M,A}\otimes A} &
             A\otimes M\otimes A \ar[r]^{A\otimes w_{M,A}} & \ar[d]^{\mu_A\otimes M} A\otimes A\otimes M\\
      M\otimes A \ar[rr]^{w_{M,A}} & & A\otimes M\\
      & \ar[lu]^{M\otimes e} M \ar[ru]_{e\otimes M}
    }
  \end{equation*}
  One can compare this diagram with the ``bow-tie'' diagram of
  entwining structures ~\cite[Diagram
  5.5]{Brzezinski:DifferentialCalculus}.  However, there one requires
  $M$ to be a coalgebra and one has similar compatibility conditions
  on the comultiplication structure.  Here we do not require $M$ to be
  a coalgebra.  The interaction between the multiplication morphism
  and $w_{M,A}$, and its inverse $w_{M,A}^{-1}$ if it exits, will be
  denoted by
\begin{equation*}
\uput{0}[0](-5.35,-0.00){M}
\uput{0}[0](-4.85,-0.00){A}
\uput{0}[0](-4.35,-0.00){A}
\uput{0}[0](-3.10,-0.00){M}
\uput{0}[0](-2.60,-0.00){A}
\uput{0}[0](-2.10,-0.00){A}
\uput{0}[0](-0.60,-0.00){A}
\uput{0}[0](-0.10,-0.00){A}
\uput{0}[0](0.40,-0.00){M}
\uput{0}[0](1.40,-0.00){A}
\uput{0}[0](1.90,-0.00){A}
\uput{0}[0](2.40,-0.00){M}
\psline(-5.2500,-0.2500)(-5.2500,-0.3750)
\psline(-4.7500,-0.3750)(-5.2500,-0.6250)
\psline[border=2pt](-5.2500,-0.3750)(-4.7500,-0.6250)
\psline(-4.7500,-0.6250)(-4.7500,-0.7500)
\psline(-4.7500,-0.2500)(-4.7500,-0.3750)
\psline(-5.2500,-0.6250)(-5.2500,-0.7500)
\psline(-4.2500,-0.2500)(-4.2500,-0.7500)
\psline(-3.0000,-0.2500)(-3.0000,-0.7500)
\psline(-2.5000,-0.2500)(-2.5000,-0.3750)
\psline(-2.5000,-0.3750)(-2.2500,-0.6250)
\psline(-2.2500,-0.6250)(-2.2500,-0.7500)
\psline(-2.0000,-0.2500)(-2.0000,-0.3750)
\psline(-2.0000,-0.3750)(-2.2500,-0.6250)
\psline(-0.5000,-0.2500)(-0.5000,-0.7500)
\psline(0.0000,-0.2500)(0.0000,-0.3750)
\psline(0.0000,-0.3750)(0.5000,-0.6250)
\psline[border=2pt](0.5000,-0.3750)(0.0000,-0.6250)
\psline(0.5000,-0.6250)(0.5000,-0.7500)
\psline(0.5000,-0.2500)(0.5000,-0.3750)
\psline(0.0000,-0.6250)(0.0000,-0.7500)
\psline(1.5000,-0.2500)(1.5000,-0.3750)
\psline(1.5000,-0.3750)(1.7500,-0.6250)
\psline(1.7500,-0.6250)(1.7500,-0.7500)
\psline(2.0000,-0.2500)(2.0000,-0.3750)
\psline(2.0000,-0.3750)(1.7500,-0.6250)
\psline(2.5000,-0.2500)(2.5000,-0.7500)
\psline(-5.2500,-0.7500)(-5.2500,-1.2500)
\psline(-4.7500,-0.7500)(-4.7500,-0.8750)
\psline(-4.2500,-0.8750)(-4.7500,-1.1250)
\psline[border=2pt](-4.7500,-0.8750)(-4.2500,-1.1250)
\psline(-4.2500,-1.1250)(-4.2500,-1.2500)
\psline(-4.2500,-0.7500)(-4.2500,-0.8750)
\psline(-4.7500,-1.1250)(-4.7500,-1.2500)
\uput{0}[0](-3.85,-1.00){~=}
\psline(-3.0000,-0.7500)(-3.0000,-0.8750)
\psline(-2.2500,-0.8750)(-3.0000,-1.1250)
\psline[border=2pt](-3.0000,-0.8750)(-2.2500,-1.1250)
\psline(-2.2500,-1.1250)(-2.2500,-1.2500)
\psline(-2.2500,-0.7500)(-2.2500,-0.8750)
\psline(-3.0000,-1.1250)(-3.0000,-1.2500)
\psline(-0.5000,-0.7500)(-0.5000,-0.8750)
\psline(-0.5000,-0.8750)(0.0000,-1.1250)
\psline[border=2pt](0.0000,-0.8750)(-0.5000,-1.1250)
\psline(0.0000,-1.1250)(0.0000,-1.2500)
\psline(0.0000,-0.7500)(0.0000,-0.8750)
\psline(-0.5000,-1.1250)(-0.5000,-1.2500)
\psline(0.5000,-0.7500)(0.5000,-1.2500)
\uput{0}[0](0.90,-1.00){~=}
\psline(1.7500,-0.7500)(1.7500,-0.8750)
\psline(1.7500,-0.8750)(2.5000,-1.1250)
\psline[border=2pt](2.5000,-0.8750)(1.7500,-1.1250)
\psline(2.5000,-1.1250)(2.5000,-1.2500)
\psline(2.5000,-0.7500)(2.5000,-0.8750)
\psline(1.7500,-1.1250)(1.7500,-1.2500)
\psline(-5.2500,-1.2500)(-5.2500,-1.3750)
\psline(-5.2500,-1.3750)(-5.0000,-1.6250)
\psline(-5.0000,-1.6250)(-5.0000,-1.7500)
\psline(-4.7500,-1.2500)(-4.7500,-1.3750)
\psline(-4.7500,-1.3750)(-5.0000,-1.6250)
\psline(-4.2500,-1.2500)(-4.2500,-1.7500)
\uput{0}[0](-3.10,-1.50){A}
\uput{0}[0](-2.35,-1.50){M}
\psline(-0.5000,-1.2500)(-0.5000,-1.7500)
\psline(0.0000,-1.2500)(0.0000,-1.3750)
\psline(0.0000,-1.3750)(0.2500,-1.6250)
\psline(0.2500,-1.6250)(0.2500,-1.7500)
\psline(0.5000,-1.2500)(0.5000,-1.3750)
\psline(0.5000,-1.3750)(0.2500,-1.6250)
\uput{0}[0](1.65,-1.50){M}
\uput{0}[0](2.40,-1.50){A}
\uput{0}[0](-5.10,-2.00){A}
\uput{0}[0](-4.35,-2.00){M}
\uput{0}[0](-0.60,-2.00){M}
\uput{0}[0](0.15,-2.00){A}
\vspace{1.9500cm}
\end{equation*}
  Similarly, the interaction between the unit morphism and $w_{M,A}$,
  and $w_{M,A}^{-1}$ if it exits, will be denoted by
\begin{equation*}
\uput{0}[0](-5.35,-0.00){M}
\uput{0}[0](-3.35,-0.00){M}
\uput{0}[0](-1.60,-0.00){M}
\uput{0}[0](-0.60,-0.00){M}
\psline(-5.2500,-0.2500)(-5.2500,-0.7500)
\psline{o}(-4.7500,-0.5000)(-4.7500,-0.7500)
\psline{o}(-3.7500,-0.5000)(-3.7500,-0.7500)
\psline(-3.2500,-0.2500)(-3.2500,-0.7500)
\psline{o}(-2.0000,-0.5000)(-2.0000,-0.7500)
\psline(-1.5000,-0.2500)(-1.5000,-0.7500)
\psline(-0.5000,-0.2500)(-0.5000,-0.7500)
\psline{o}(0.0000,-0.5000)(0.0000,-0.7500)
\psline(-5.2500,-0.7500)(-5.2500,-0.8750)
\psline(-4.7500,-0.8750)(-5.2500,-1.1250)
\psline[border=2pt](-5.2500,-0.8750)(-4.7500,-1.1250)
\psline(-4.7500,-1.1250)(-4.7500,-1.2500)
\psline(-4.7500,-0.7500)(-4.7500,-0.8750)
\psline(-5.2500,-1.1250)(-5.2500,-1.2500)
\uput{0}[0](-4.35,-1.00){=}
\uput{0}[0](-3.85,-1.00){A}
\uput{0}[0](-3.35,-1.00){M}
\psline(-2.0000,-0.7500)(-2.0000,-0.8750)
\psline(-2.0000,-0.8750)(-1.5000,-1.1250)
\psline[border=2pt](-1.5000,-0.8750)(-2.0000,-1.1250)
\psline(-1.5000,-1.1250)(-1.5000,-1.2500)
\psline(-1.5000,-0.7500)(-1.5000,-0.8750)
\psline(-2.0000,-1.1250)(-2.0000,-1.2500)
\uput{0}[0](-1.10,-1.00){=}
\uput{0}[0](-0.60,-1.00){M}
\uput{0}[0](-0.10,-1.00){A}
\uput{0}[0](-5.35,-1.50){A}
\uput{0}[0](-4.85,-1.50){M}
\uput{0}[0](-2.10,-1.50){M}
\uput{0}[0](-1.60,-1.50){A}
\vspace{1.4500cm}
\end{equation*}
  A $w$-transpositive coalgebra $(C,\delta_C,\varepsilon)$ is a
  $w^{op}$-transpositive algebra in the opposite monoidal category
  $(\C{C}^{op},\otimes^{op})$ where $U\otimes^{op}V\assign V\otimes U$
  for any two objects $U,V\in Ob(\C{C}^{op})$.
\end{definition}

For the examples we are going to consider below, we fix a commutative
associative unital ring $k$.  Our base symmetric monoidal category is
$\lmod{k}$ the category of $k$-modules with $\otimes_k$ the ordinary
tensor product over $k$ taken as the monoidal product $\otimes$.  We
also assume $B$ is an associative/coassociative unital/counital
bialgebra, or a Hopf algebra with an invertible antipode whenever it
is necessary.

\begin{example}
  Fix a right/left $B$-module/comodule $(M,\beta_M,\lambda_M)$.  We
  let $w_{M,X}:M\otimes X\xra{}X\otimes M$ be a transposition if (i)
  $(X,\alpha_X)$ is a left $B$-module and (ii) $w_{M,X}$ is defined by
  the formula
  \[ w_{M,X}(m\otimes x) \assign m_{(-1)}x\otimes m_{(0)} 
     = (\alpha_X\otimes M)\circ(B\otimes s_{M,X})\circ(\lambda_M\otimes X)(m\otimes x)     
  \] 
  for any $m\otimes x\in M\otimes X$ where $s_{M,X}$ is the ordinary
  switch morphism.  For this class of transpositions $w$, an algebra
  $(A,\mu_A,e)$ is $w$-transpositive if $A$ is a left $B$-module algebra.
  Similarly $(C,\delta_C,\varepsilon)$ is a $w$-transpositive coalgebra, if
  $C$ is a right $B$-comodule coalgebra.
\end{example}

\begin{example}
  Fix a left/right $B$-module/comodule $(M,\alpha_M,\rho_M)$.  We let
  $w_{M,X}:M\otimes X\xra{}X\otimes M$ be a transposition if (i)
  $(X,\beta_X)$ is a right $B$-comodule and (ii) $w_{M,X}$ is defined
  by the formula
  \[ w_{M,X}(m\otimes x) \assign x_{(0)}\otimes m x_{(1)} =
  (\alpha_X\otimes M)\circ(s_{M,X}\otimes B)\circ(M\otimes\beta_X)(m\otimes x)
  \] 
  for any $m\otimes x\in M\otimes X$ where $s_{M,X}$ is the ordinary
  switch morphism.  For this class of transpositions $w$, an algebra
  $(A,\mu_A,e)$ is $w$-transpositive if $A$ is a right $B$-comodule
  algebra.  Similarly $(C,\delta_C,\varepsilon)$ is a $w$-transpositive
  coalgebra, if $C$ is a left $B$-module coalgebra.
\end{example}

\section{The universal para-(co)cyclic theory} 
\label{sec:UniversalParaCocyclicTheory}

\begin{definition}
  Let $S$ be the category with objects $\{0,1\}$ where there is one
  unique morphism $i\xra{}j$ between any two objects $i,j\in\{0,1\}$.
  A functor $F\colon S\xra{}\C{C}$ will be called an $S$-module.
\end{definition}

\begin{lemma}
  \label{lem:UniqueLift}Let $\C{G}$ and $\C{G}'$ be two
  groupoids which have the property that between any two objects there
  is a unique morphism.  Let $F,G \colon \C{G} \longrightarrow \C{C}$
  and $F' \colon \C{G}' \longrightarrow \C{C}$ be three arbitrary
  functors and let $g, h \in {\rm Ob} ( \C{G} )$ and $g' \in {\rm Ob}
  ( \C{G'} )$ be three arbitrary objects.
  \begin{enumerate}
  \item Any morphism $F ( g ) \xra{u} G ( h )$ in $\C{C}$ can be
    lifted to a natural transformation of functors of the form
    $F\xra{u} G$.
    
  \item Any morphism $F ( g ) \xra{v} F' ( g' )$ in $\C{C}$ can be
    lifted to a morphism in $\C{C}$ of the form $\colim_{\C{G}} F
    \xra{v} \colim_{\C{G}'} F'$
  \end{enumerate}
\end{lemma}

\begin{proof}
  We will denote the unique morphism between from an object $x$ to
  another object $y$ in $\C{G}$ by $x\xra{\alpha_{y,x}}y$.  Let
  $F(g)\xra{u}G(h)$ be an arbitrary morphism in $\C{C}$ and define
  \[ u_x \assign G(\alpha_{x,h})\circ u\circ F(\alpha_{g,x}) \]
  for any $x\in Ob(\C{G})$.  In order $u_\cdot$ to define a natural
  transformation, for any $x,y\in Ob(\C{G})$ one must have
  \[ u_y \circ F(\alpha_{y,x}) = G(\alpha_{y,x})\circ u_x \]
  So we check
  \begin{align*}
    u_y\circ F(\alpha_{y,x}) 
    = & G(\alpha_{y,h})\circ u\circ F(\alpha_{g,y})\circ F(\alpha_{y,x})\\
    = & G(\alpha_{y,x})\circ G(\alpha_{x,h})\circ u\circ F(\alpha_{g,x})
    = G(\alpha_{y,x})\circ u_x
  \end{align*}
  for any $x,y\in Ob(\C{G})$ as we wanted to show.  This finishes the
  first part of the assertion.

  For the second part, let $F(x)\xra{\phi(x)}\colim_\C{G}F$ and
  $F'(x')\xra{\phi'(x')}\colim_{\C{G}'}F'$ be the structure
  morphisms of the corresponding colimits.  Then one has morphisms 
  \[ \psi(x) \assign \phi'(g')\circ v \circ F(\alpha_{g,x}) \]
  of the form $x\xra{\psi(x)}\colim_{\C{G}'}F'$ for any $x\in
  Ob(\C{G})$.  We check that
  \[ \psi(x)\circ F(\alpha_{x,y}) =
     \phi'(g')\circ v \circ F(\alpha_{g,x})\circ F(\alpha_{x,y})
     = \phi'(g')\circ v \circ F(\alpha_{g,y}) = \psi(y)
  \]
  for any $x,y\in Ob(\C{G})$ meaning there is a unique morphism
  $\colim_\C{G}F\xra{}\colim_{\C{G}'}G$ which, by abuse of notation,
  we still denote by $v$.
\end{proof}

\begin{definition}
  Let $C$ and $M$ be two arbitrary objects in $\C{C}$ such that we
  have a transposition $M\otimes C\xra{w_{M,C}}C\otimes M$.  For every
  $n\geq 0$, we define an $S$-module $P_n ( C, M )$ in $\C{C}$ as
  follows: let $P_n( C, M )$ is the functor from $S$ to $\C{C}$
  given on the objects by
  \begin{align}
    P_n (C,M) (0) \assign & C^{\otimes n} \otimes M\otimes C &
    P_n (C,M) (1) \assign & C^{\otimes n+1} \otimes M
  \end{align}
  Moreover, 
  \begin{equation}
    P_n(C,M)(0\xra{}1) \assign 
    \left(C^{\otimes n} \otimes M\otimes C\xra{t_{n+2}} 
          C^{\otimes n+1} \otimes M\right)
  \end{equation}
  is the cyclic permutation coming from the symmetric monoidal
  structure of $\C{C}$ thus its inverse provides
  $P_n(C,M)(1\xra{}0)$.  
\end{definition}

\begin{definition}
  Let $\Lambda$ be Connes' cyclic category {\cite{Connes:ExtFunctors}}
  and $\Lambda_{\B{N}}$ and $\Lambda_{\B{Z}}$ be the variations of
  $\Lambda$ as defined in {\cite{Kaygun:BivariantHopf}}.  Let us
  recall the presentation we will use in this paper: the category
  $\Lambda_\B{N}$ has objects $[n]$ indexed by natural numbers $n\geq
  0$ and is generated by morphisms $[n]\xra{\partial^n_j}[n+1]$,
  $[n+1]\xra{\sigma^n_i}[n]$ and $[n]\xra{\tau_n^\ell}[n]$ with $0\leq
  j\leq n+1$, $0\leq i\leq n$ and $\ell\in\B{N}$.  These generators
  are subject to the following relations
  \begin{align*}
    \partial^{n+1}_i\partial^n_j = & \partial^{n+1}_{j+1}\partial^n_i \text{ \ \ and \ \ }
    \sigma^{n-1}_j\sigma^n_i = \sigma^{n-1}_i\sigma^n_{j+1} \text{ for } i\leq j \ \  \text{ and }
    \tau_n^s\tau_n^t = \tau_n^{s+t} \text{ for } s,t\in\B{N}\\
    \sigma^n_i\partial^n_i = & \sigma^n_i\partial^n_{i+1} = \tau_n^0 \text{ \ \ and\ \ }
    \partial^n_i\sigma^n_j 
      = \begin{cases}
	\sigma^{n+1}_{j+1}\partial^{n+1}_i & \text{ if }i\leq j\\
	\sigma^{n+1}_j\partial^{n+1}_{i+1} & \text{ if }i>j
	\end{cases}\\
    \partial^n_j\tau_n^i = & \tau_{n+1}^{i+p}\partial^n_q 
        \text{ where $(i+j)=(n+1)p+q$ with $0\leq q\leq n$}\\
    \tau_n^i\sigma^n_j = & \sigma^n_q\tau_{n+1}^{i+p} 
        \text{ where $(-i+j)=(n+1)(-p)+q$ with $0\leq q\leq n$}
  \end{align*}
  The category $\Lambda_\B{Z}$ is an extension of $\Lambda_\B{N}$
  where we allow morphisms of the form $\tau_n^i$ with $i\in\B{Z}$.
  Then $\Lambda$ is a quotient of $\Lambda_\B{Z}$ where we put the
  extra relations $\tau_n^{n+1}=id_n$ for $n\geq 0$.  The category
  $\Lambda_+$ is the subcategory of $\Lambda_\B{N}$ generated by
  $\partial^n_j$ and $\sigma^n_i$ with only $0\leq i\leq n$ and $0\leq
  j\leq n$.  A functor $F \colon \Lambda \longrightarrow \C{C}$ will
  be referred as {\em a cocyclic module in} $\C{C}$ while any functor
  of the form $F \colon \Lambda_{\B{N}} \longrightarrow \C{C}$ or $F
  \colon \Lambda_{\B{Z}} \longrightarrow \C{C}$ will be referred as
  {\em a para-cocyclic module in} $\C{C}$.  A (para-)cyclic module $F$
  in $\C{C}$ is defined to be a (para-)cocyclic module in
  $\C{C}^{op}$.  A morphism between (para-)cocyclic modules $h \colon
  F \longrightarrow G$ in $\C{C}$ is just a natural transformation of
  functors.
\end{definition}

\begin{theorem}\label{thm:HomotopyParaCyclic}
  Let $(C,\delta_C,\varepsilon)$ be a $w$-transpositive coalgebra.
  Let $\colim_S P_\bullet(C,M)$ be level-wise colimit of
  $P_\bullet(C,M)$.  Then $\colim_S P_\bullet ( C, M )$
  carries a para-cocyclic module structure.
\end{theorem}

\begin{proof}
  The cosimplicial structure morphisms are given by
  \[ \partial_i \assign C^{\otimes i} \otimes \delta_C \otimes C^{\otimes n - i} 
         \otimes M 
     {\rm \ \ \ and\ \ \ } 
     \sigma_j \assign C^{\otimes j + 1} \otimes \varepsilon \otimes 
         C^{\otimes n - 1 - j} \otimes M 
  \]
  which are defined only for $0 \leqslant i \leqslant n$ and $0
  \leqslant j \leqslant n - 1$ and on $C^{\otimes n + 1} \otimes M$.
  We also let
  \[ \partial_{n + 1} \assign ( C^{\otimes n} \otimes w_{M, C} \otimes C )
     \circ ( C^{\otimes n} \otimes M \otimes \delta_C ) 
  \]
  which is a morphism defined on $C^{\otimes n}\otimes M\otimes C$.
  The fact that the morphisms $\partial_i$ and $\sigma_j$ for $0
  \leqslant i \leqslant n + 1$ and $0 \leqslant j \leqslant n$ are
  well-defined on the level-wise colimits follows from Lemma
  \ref{lem:UniqueLift}.  The para-cocyclic structure morphisms are
  already used in this definition since we are going to define
  \[ \tau_n \assign C^{\otimes n} \otimes w_{M, C} \] 
  for any $n \geqslant 0$.  This is a morphism of the form
  \[ P_n(C,M)(0)\xra{\tau_n}P_n(C,M)(1) \]
  The fact that $\tau_n$ is well-defined on $\colim_S P_n (A,M)$ for
  any $n\geq 0$ again is a consequence of Lemma \ref{lem:UniqueLift}.
  The verification of the cosimplicial identities between $\partial_i$
  and $\sigma_j$ for the range $0 \leqslant i \leqslant n$ and $0
  \leqslant j \leqslant n$ is standard and follows from the fact that
  $C$ is a coassociative counital coalgebra in $\C{C}$.  Next, we
  consider $\partial_j\partial_{n+1}$.  If $0\leq j\leq n$, one can
  describe the composition by
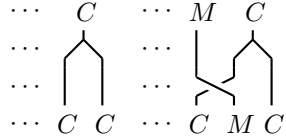
\begin{equation*}
\uput{0}[0](-5.35,-0.00){~\cdots}
\uput{0}[0](-4.35,-0.00){C}
\uput{0}[0](-3.60,-0.00){~\cdots}
\uput{0}[0](-2.85,-0.00){M}
\uput{0}[0](-2.10,-0.00){C}
\uput{0}[0](-5.35,-0.50){~\cdots}
\psline(-4.2500,-0.2500)(-4.2500,-0.3750)
\psline(-4.2500,-0.3750)(-4.5000,-0.6250)
\psline(-4.5000,-0.6250)(-4.5000,-0.7500)
\psline(-4.2500,-0.3750)(-4.0000,-0.6250)
\psline(-4.0000,-0.6250)(-4.0000,-0.7500)
\uput{0}[0](-3.60,-0.50){~\cdots}
\psline(-2.7500,-0.2500)(-2.7500,-0.7500)
\psline(-2.0000,-0.2500)(-2.0000,-0.3750)
\psline(-2.0000,-0.3750)(-2.2500,-0.6250)
\psline(-2.2500,-0.6250)(-2.2500,-0.7500)
\psline(-2.0000,-0.3750)(-1.7500,-0.6250)
\psline(-1.7500,-0.6250)(-1.7500,-0.7500)
\uput{0}[0](-5.35,-1.00){~\cdots}
\psline(-4.5000,-0.7500)(-4.5000,-1.2500)
\psline(-4.0000,-0.7500)(-4.0000,-1.2500)
\uput{0}[0](-3.60,-1.00){~\cdots}
\psline(-2.7500,-0.7500)(-2.7500,-0.8750)
\psline(-2.2500,-0.8750)(-2.7500,-1.1250)
\psline[border=2pt](-2.7500,-0.8750)(-2.2500,-1.1250)
\psline(-2.2500,-1.1250)(-2.2500,-1.2500)
\psline(-2.2500,-0.7500)(-2.2500,-0.8750)
\psline(-2.7500,-1.1250)(-2.7500,-1.2500)
\psline(-1.7500,-0.7500)(-1.7500,-1.2500)
\uput{0}[0](-5.35,-1.50){~\cdots}
\uput{0}[0](-4.60,-1.50){C}
\uput{0}[0](-4.10,-1.50){C}
\uput{0}[0](-3.60,-1.50){~\cdots}
\uput{0}[0](-2.85,-1.50){C}
\uput{0}[0](-2.35,-1.50){M}
\uput{0}[0](-1.85,-1.50){C}
\vspace{1.4500cm}
\end{equation*}
  This shows $\partial_j\partial_{n+1} = \partial_{n+2}\partial_j$ for
  $0\leq j\leq n$.  For $j=n+1$, by using the fact that $C$ is a
  $w$-transpositive coalgebra we see that
  $\partial_{n+1}\partial_{n+1}$ can be described as
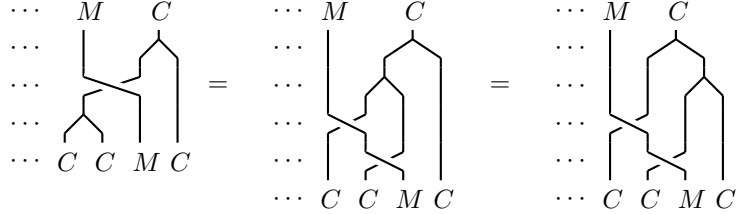
\begin{equation*}
\uput{0}[0](-5.35,-0.00){~\cdots}
\uput{0}[0](-4.35,-0.00){M}
\uput{0}[0](-3.35,-0.00){C}
\uput{0}[0](-1.85,-0.00){~\cdots}
\uput{0}[0](-1.10,-0.00){M}
\uput{0}[0](0.02,-0.00){C}
\uput{0}[0](1.90,-0.00){~\cdots}
\uput{0}[0](2.65,-0.00){M}
\uput{0}[0](3.52,-0.00){C}
\uput{0}[0](-5.35,-0.50){~\cdots}
\psline(-4.2500,-0.2500)(-4.2500,-0.7500)
\psline(-3.2500,-0.2500)(-3.2500,-0.3750)
\psline(-3.2500,-0.3750)(-3.5000,-0.6250)
\psline(-3.5000,-0.6250)(-3.5000,-0.7500)
\psline(-3.2500,-0.3750)(-3.0000,-0.6250)
\psline(-3.0000,-0.6250)(-3.0000,-0.7500)
\uput{0}[0](-1.85,-0.50){~\cdots}
\psline(-1.0000,-0.2500)(-1.0000,-0.7500)
\psline(0.1250,-0.2500)(0.1250,-0.3750)
\psline(0.1250,-0.3750)(-0.2500,-0.6250)
\psline(-0.2500,-0.6250)(-0.2500,-0.7500)
\psline(0.1250,-0.3750)(0.5000,-0.6250)
\psline(0.5000,-0.6250)(0.5000,-0.7500)
\uput{0}[0](1.90,-0.50){~\cdots}
\psline(2.7500,-0.2500)(2.7500,-0.7500)
\psline(3.6250,-0.2500)(3.6250,-0.3750)
\psline(3.6250,-0.3750)(3.2500,-0.6250)
\psline(3.2500,-0.6250)(3.2500,-0.7500)
\psline(3.6250,-0.3750)(4.0000,-0.6250)
\psline(4.0000,-0.6250)(4.0000,-0.7500)
\uput{0}[0](-5.35,-1.00){~\cdots}
\psline(-4.2500,-0.7500)(-4.2500,-0.8750)
\psline(-3.5000,-0.8750)(-4.2500,-1.1250)
\psline[border=2pt](-4.2500,-0.8750)(-3.5000,-1.1250)
\psline(-3.5000,-1.1250)(-3.5000,-1.2500)
\psline(-3.5000,-0.7500)(-3.5000,-0.8750)
\psline(-4.2500,-1.1250)(-4.2500,-1.2500)
\psline(-3.0000,-0.7500)(-3.0000,-1.2500)
\uput{0}[0](-2.60,-1.00){=}
\uput{0}[0](-1.85,-1.00){~\cdots}
\psline(-1.0000,-0.7500)(-1.0000,-1.2500)
\psline(-0.2500,-0.7500)(-0.2500,-0.8750)
\psline(-0.2500,-0.8750)(-0.5000,-1.1250)
\psline(-0.5000,-1.1250)(-0.5000,-1.2500)
\psline(-0.2500,-0.8750)(0.0000,-1.1250)
\psline(0.0000,-1.1250)(0.0000,-1.2500)
\psline(0.5000,-0.7500)(0.5000,-1.2500)
\uput{0}[0](1.15,-1.00){=}
\uput{0}[0](1.90,-1.00){~\cdots}
\psline(2.7500,-0.7500)(2.7500,-1.2500)
\psline(3.2500,-0.7500)(3.2500,-1.2500)
\psline(4.0000,-0.7500)(4.0000,-0.8750)
\psline(4.0000,-0.8750)(3.7500,-1.1250)
\psline(3.7500,-1.1250)(3.7500,-1.2500)
\psline(4.0000,-0.8750)(4.2500,-1.1250)
\psline(4.2500,-1.1250)(4.2500,-1.2500)
\uput{0}[0](-5.35,-1.50){~\cdots}
\psline(-4.2500,-1.2500)(-4.2500,-1.3750)
\psline(-4.2500,-1.3750)(-4.5000,-1.6250)
\psline(-4.5000,-1.6250)(-4.5000,-1.7500)
\psline(-4.2500,-1.3750)(-4.0000,-1.6250)
\psline(-4.0000,-1.6250)(-4.0000,-1.7500)
\psline(-3.5000,-1.2500)(-3.5000,-1.7500)
\psline(-3.0000,-1.2500)(-3.0000,-1.7500)
\uput{0}[0](-1.85,-1.50){~\cdots}
\psline(-1.0000,-1.2500)(-1.0000,-1.3750)
\psline(-0.5000,-1.3750)(-1.0000,-1.6250)
\psline[border=2pt](-1.0000,-1.3750)(-0.5000,-1.6250)
\psline(-0.5000,-1.6250)(-0.5000,-1.7500)
\psline(-0.5000,-1.2500)(-0.5000,-1.3750)
\psline(-1.0000,-1.6250)(-1.0000,-1.7500)
\psline(0.0000,-1.2500)(0.0000,-1.7500)
\psline(0.5000,-1.2500)(0.5000,-1.7500)
\uput{0}[0](1.90,-1.50){~\cdots}
\psline(2.7500,-1.2500)(2.7500,-1.3750)
\psline(3.2500,-1.3750)(2.7500,-1.6250)
\psline[border=2pt](2.7500,-1.3750)(3.2500,-1.6250)
\psline(3.2500,-1.6250)(3.2500,-1.7500)
\psline(3.2500,-1.2500)(3.2500,-1.3750)
\psline(2.7500,-1.6250)(2.7500,-1.7500)
\psline(3.7500,-1.2500)(3.7500,-1.7500)
\psline(4.2500,-1.2500)(4.2500,-1.7500)
\uput{0}[0](-5.35,-2.00){~\cdots}
\uput{0}[0](-4.60,-2.00){C}
\uput{0}[0](-4.10,-2.00){C}
\uput{0}[0](-3.60,-2.00){M}
\uput{0}[0](-3.10,-2.00){C}
\uput{0}[0](-1.85,-2.00){~\cdots}
\psline(-1.0000,-1.7500)(-1.0000,-2.2500)
\psline(-0.5000,-1.7500)(-0.5000,-1.8750)
\psline(0.0000,-1.8750)(-0.5000,-2.1250)
\psline[border=2pt](-0.5000,-1.8750)(0.0000,-2.1250)
\psline(0.0000,-2.1250)(0.0000,-2.2500)
\psline(0.0000,-1.7500)(0.0000,-1.8750)
\psline(-0.5000,-2.1250)(-0.5000,-2.2500)
\psline(0.5000,-1.7500)(0.5000,-2.2500)
\uput{0}[0](1.90,-2.00){~\cdots}
\psline(2.7500,-1.7500)(2.7500,-2.2500)
\psline(3.2500,-1.7500)(3.2500,-1.8750)
\psline(3.7500,-1.8750)(3.2500,-2.1250)
\psline[border=2pt](3.2500,-1.8750)(3.7500,-2.1250)
\psline(3.7500,-2.1250)(3.7500,-2.2500)
\psline(3.7500,-1.7500)(3.7500,-1.8750)
\psline(3.2500,-2.1250)(3.2500,-2.2500)
\psline(4.2500,-1.7500)(4.2500,-2.2500)
\uput{0}[0](-1.85,-2.50){~\cdots}
\uput{0}[0](-1.10,-2.50){C}
\uput{0}[0](-0.60,-2.50){C}
\uput{0}[0](-0.10,-2.50){M}
\uput{0}[0](0.40,-2.50){C}
\uput{0}[0](1.90,-2.50){~\cdots}
\uput{0}[0](2.65,-2.50){C}
\uput{0}[0](3.15,-2.50){C}
\uput{0}[0](3.65,-2.50){M}
\uput{0}[0](4.15,-2.50){C}
\vspace{2.4500cm}
\end{equation*}
 which is equivalent to saying $\partial_{n+1}\partial_{n+1} =
 \partial_{n+2}\partial_{n+1}$.  This finishes the proof that
 $\colim_S P_\bullet (C,M)$ is pre-cosimplicial.  Now we consider
 $\sigma_i\partial_{n+1}$.  If $0\leq i< n$, the composition can be
 described by
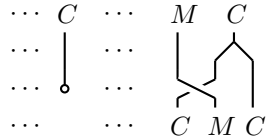
\begin{equation*}
\uput{0}[0](-5.35,-0.00){~\cdots}
\uput{0}[0](-4.60,-0.00){C}
\uput{0}[0](-4.10,-0.00){~\cdots}
\uput{0}[0](-3.10,-0.00){M}
\uput{0}[0](-2.35,-0.00){C}
\uput{0}[0](-5.35,-0.50){~\cdots}
\psline(-4.5000,-0.2500)(-4.5000,-0.7500)
\uput{0}[0](-4.10,-0.50){~\cdots}
\psline(-3.0000,-0.2500)(-3.0000,-0.7500)
\psline(-2.2500,-0.2500)(-2.2500,-0.3750)
\psline(-2.2500,-0.3750)(-2.5000,-0.6250)
\psline(-2.5000,-0.6250)(-2.5000,-0.7500)
\psline(-2.2500,-0.3750)(-2.0000,-0.6250)
\psline(-2.0000,-0.6250)(-2.0000,-0.7500)
\uput{0}[0](-5.35,-1.00){~\cdots}
\psline{-o}(-4.5000,-0.7500)(-4.5000,-1.0000)
\uput{0}[0](-4.10,-1.00){~\cdots}
\psline(-3.0000,-0.7500)(-3.0000,-0.8750)
\psline(-2.5000,-0.8750)(-3.0000,-1.1250)
\psline[border=2pt](-3.0000,-0.8750)(-2.5000,-1.1250)
\psline(-2.5000,-1.1250)(-2.5000,-1.2500)
\psline(-2.5000,-0.7500)(-2.5000,-0.8750)
\psline(-3.0000,-1.1250)(-3.0000,-1.2500)
\psline(-2.0000,-0.7500)(-2.0000,-1.2500)
\uput{0}[0](-5.35,-1.50){~\cdots}
\uput{0}[0](-4.10,-1.50){~\cdots}
\uput{0}[0](-3.10,-1.50){C}
\uput{0}[0](-2.60,-1.50){M}
\uput{0}[0](-2.10,-1.50){C}
\vspace{1.4500cm}
\end{equation*}
 Then one can easily see that $\sigma_i\partial_{n+1} =
 \partial_n\sigma_i$ for $0\leq i< n$.  We also observe that
 $\sigma_n\partial_{n+1} = id$ since
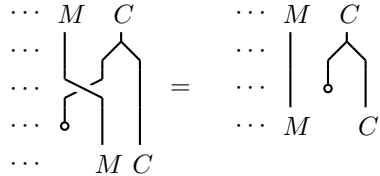
\begin{equation*}
\uput{0}[0](-5.35,-0.00){~\cdots}
\uput{0}[0](-4.60,-0.00){M}
\uput{0}[0](-3.85,-0.00){C}
\uput{0}[0](-2.35,-0.00){~\cdots}
\uput{0}[0](-1.60,-0.00){M}
\uput{0}[0](-0.85,-0.00){C}
\uput{0}[0](-5.35,-0.50){~\cdots}
\psline(-4.5000,-0.2500)(-4.5000,-0.7500)
\psline(-3.7500,-0.2500)(-3.7500,-0.3750)
\psline(-3.7500,-0.3750)(-4.0000,-0.6250)
\psline(-4.0000,-0.6250)(-4.0000,-0.7500)
\psline(-3.7500,-0.3750)(-3.5000,-0.6250)
\psline(-3.5000,-0.6250)(-3.5000,-0.7500)
\uput{0}[0](-2.35,-0.50){~\cdots}
\psline(-1.5000,-0.2500)(-1.5000,-0.7500)
\psline(-0.7500,-0.2500)(-0.7500,-0.3750)
\psline(-0.7500,-0.3750)(-1.0000,-0.6250)
\psline(-1.0000,-0.6250)(-1.0000,-0.7500)
\psline(-0.7500,-0.3750)(-0.5000,-0.6250)
\psline(-0.5000,-0.6250)(-0.5000,-0.7500)
\uput{0}[0](-5.35,-1.00){~\cdots}
\psline(-4.5000,-0.7500)(-4.5000,-0.8750)
\psline(-4.0000,-0.8750)(-4.5000,-1.1250)
\psline[border=2pt](-4.5000,-0.8750)(-4.0000,-1.1250)
\psline(-4.0000,-1.1250)(-4.0000,-1.2500)
\psline(-4.0000,-0.7500)(-4.0000,-0.8750)
\psline(-4.5000,-1.1250)(-4.5000,-1.2500)
\psline(-3.5000,-0.7500)(-3.5000,-1.2500)
\uput{0}[0](-3.10,-1.00){=}
\uput{0}[0](-2.35,-1.00){~\cdots}
\psline(-1.5000,-0.7500)(-1.5000,-1.2500)
\psline{-o}(-1.0000,-0.7500)(-1.0000,-1.0000)
\psline(-0.5000,-0.7500)(-0.5000,-1.2500)
\uput{0}[0](-5.35,-1.50){~\cdots}
\psline{-o}(-4.5000,-1.2500)(-4.5000,-1.5000)
\psline(-4.0000,-1.2500)(-4.0000,-1.7500)
\psline(-3.5000,-1.2500)(-3.5000,-1.7500)
\uput{0}[0](-2.35,-1.50){~\cdots}
\uput{0}[0](-1.60,-1.50){M}
\uput{0}[0](-0.60,-1.50){C}
\uput{0}[0](-5.35,-2.00){~\cdots}
\uput{0}[0](-4.10,-2.00){M}
\uput{0}[0](-3.60,-2.00){C}
\vspace{1.9500cm}
\end{equation*}
 This finishes the proof that $\colim_S P_\bullet(C,M)$ is a
 cosimplicial object in $\C{C}$.  Now we must check the para-cocyclic
 identities.  First we observe that $\tau_{n+1}\partial_0 =
 \partial_{n+1}$ by definition.  Next, we consider
 $\tau_{n+1}\partial_i$.  For the range $0< i< n$ we represent the
 composition by
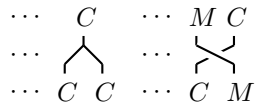
\begin{equation*}
\uput{0}[0](-5.35,-0.00){~\cdots}
\uput{0}[0](-4.35,-0.00){C}
\uput{0}[0](-3.60,-0.00){~\cdots}
\uput{0}[0](-2.85,-0.00){M}
\uput{0}[0](-2.35,-0.00){C}
\uput{0}[0](-5.35,-0.50){~\cdots}
\psline(-4.2500,-0.2500)(-4.2500,-0.3750)
\psline(-4.2500,-0.3750)(-4.5000,-0.6250)
\psline(-4.5000,-0.6250)(-4.5000,-0.7500)
\psline(-4.2500,-0.3750)(-4.0000,-0.6250)
\psline(-4.0000,-0.6250)(-4.0000,-0.7500)
\uput{0}[0](-3.60,-0.50){~\cdots}
\psline(-2.7500,-0.2500)(-2.7500,-0.3750)
\psline(-2.2500,-0.3750)(-2.7500,-0.6250)
\psline[border=2pt](-2.7500,-0.3750)(-2.2500,-0.6250)
\psline(-2.2500,-0.6250)(-2.2500,-0.7500)
\psline(-2.2500,-0.2500)(-2.2500,-0.3750)
\psline(-2.7500,-0.6250)(-2.7500,-0.7500)
\uput{0}[0](-5.35,-1.00){~\cdots}
\uput{0}[0](-4.60,-1.00){C}
\uput{0}[0](-4.10,-1.00){C}
\uput{0}[0](-3.60,-1.00){~\cdots}
\uput{0}[0](-2.85,-1.00){C}
\uput{0}[0](-2.35,-1.00){M}
\vspace{0.9500cm}
\end{equation*}
 which means one has $\tau_{n+1}\partial_i = \partial_{i-1}\tau_n$ for
 the range $0< i< n$.  For $i=n$ we consider $\partial_n\tau_n$
 which is represented by
\begin{equation*}
\uput{0}[0](-5.35,-0.00){~\cdots}
\uput{0}[0](-4.35,-0.00){M}
\uput{0}[0](-3.60,-0.00){C}
\uput{0}[0](-2.10,-0.00){~\cdots}
\uput{0}[0](-1.35,-0.00){M}
\uput{0}[0](-0.60,-0.00){C}
\uput{0}[0](-5.35,-0.50){~\cdots}
\psline(-4.2500,-0.2500)(-4.2500,-0.7500)
\psline(-3.5000,-0.2500)(-3.5000,-0.7500)
\uput{0}[0](-2.10,-0.50){~\cdots}
\psline(-1.2500,-0.2500)(-1.2500,-0.7500)
\psline(-0.5000,-0.2500)(-0.5000,-0.3750)
\psline(-0.5000,-0.3750)(-0.7500,-0.6250)
\psline(-0.7500,-0.6250)(-0.7500,-0.7500)
\psline(-0.5000,-0.3750)(-0.2500,-0.6250)
\psline(-0.2500,-0.6250)(-0.2500,-0.7500)
\uput{0}[0](-5.35,-1.00){~\cdots}
\psline(-4.2500,-0.7500)(-4.2500,-0.8750)
\psline(-3.5000,-0.8750)(-4.2500,-1.1250)
\psline[border=2pt](-4.2500,-0.8750)(-3.5000,-1.1250)
\psline(-3.5000,-1.1250)(-3.5000,-1.2500)
\psline(-3.5000,-0.7500)(-3.5000,-0.8750)
\psline(-4.2500,-1.1250)(-4.2500,-1.2500)
\uput{0}[0](-2.85,-1.00){=}
\uput{0}[0](-2.10,-1.00){~\cdots}
\psline(-1.2500,-0.7500)(-1.2500,-0.8750)
\psline(-0.7500,-0.8750)(-1.2500,-1.1250)
\psline[border=2pt](-1.2500,-0.8750)(-0.7500,-1.1250)
\psline(-0.7500,-1.1250)(-0.7500,-1.2500)
\psline(-0.7500,-0.7500)(-0.7500,-0.8750)
\psline(-1.2500,-1.1250)(-1.2500,-1.2500)
\psline(-0.2500,-0.7500)(-0.2500,-1.2500)
\uput{0}[0](-5.35,-1.50){~\cdots}
\psline(-4.2500,-1.2500)(-4.2500,-1.3750)
\psline(-4.2500,-1.3750)(-4.5000,-1.6250)
\psline(-4.5000,-1.6250)(-4.5000,-1.7500)
\psline(-4.2500,-1.3750)(-4.0000,-1.6250)
\psline(-4.0000,-1.6250)(-4.0000,-1.7500)
\psline(-3.5000,-1.2500)(-3.5000,-1.7500)
\uput{0}[0](-2.10,-1.50){~\cdots}
\psline(-1.2500,-1.2500)(-1.2500,-1.7500)
\psline(-0.7500,-1.2500)(-0.7500,-1.3750)
\psline(-0.2500,-1.3750)(-0.7500,-1.6250)
\psline[border=2pt](-0.7500,-1.3750)(-0.2500,-1.6250)
\psline(-0.2500,-1.6250)(-0.2500,-1.7500)
\psline(-0.2500,-1.2500)(-0.2500,-1.3750)
\psline(-0.7500,-1.6250)(-0.7500,-1.7500)
\uput{0}[0](-5.35,-2.00){~\cdots}
\uput{0}[0](-4.60,-2.00){C}
\uput{0}[0](-4.10,-2.00){C}
\uput{0}[0](-3.60,-2.00){M}
\uput{0}[0](-2.10,-2.00){~\cdots}
\uput{0}[0](-1.35,-2.00){C}
\uput{0}[0](-0.85,-2.00){C}
\uput{0}[0](-0.35,-2.00){M}
\vspace{1.9500cm}
\end{equation*}
 which is equivalent to saying $\partial_n\tau_n =
 \tau_{n+1}^2\partial_0 = \tau_{n+1}\partial_{n+1}$.  So far we have
 the following relations
 \[ \partial_i\tau_n = \tau_{n+1}\partial_{i+1} \text{ for } 0\leq i< n 
    \text{ and }\partial_n\tau_n = \tau_{n+1}^2\partial_0
 \]
 Using these relations one can show that
 \[ \partial_i \tau_n^j = \tau_{n+1}^{j+p} \partial_q \text{ where } (i+j) = (n+1)p + q 
 \] 
 for any $n\geq 0$, $0\leq i\leq n+1$ and $j\geq 0$ i.e.  $\colim_S
 P_\bullet(C,M)$ is a pre-para-cocyclic object in $\C{C}$.  We leave
 the verification of the identities
 \[ \tau_n^j\sigma_i = \sigma_q\tau_{n+1}^{i+p} \text{ where } (i-j) = (n+1)(-p) + q
 \] involving para-cyclic operators and the codegeneracy operators to
 the reader.
\end{proof}

For simplicity, the para-cocyclic module $\colim_S P_\bullet(C,M)$
will be denoted by $T_\bullet(C,M)$.

\section{Approximation theorems for pseudo-para-(co)cyclic objects}
\label{sec:Approximations}

\begin{definition}\label{defn:Approximation}
  Here we assume $\C{C}$ is an arbitrary small category and let
  $\C{D}$ be a subcategory.  For an arbitrary object $X$ of $\C{C}$ a
  morphism $\app{\C{D}}{X}\xra{u_X}X$ is called the approximation of
  $X$ within $\C{D}$ if (i) $\app{\C{D}}{X}$ is an object in $\C{D}$
  and (ii) every morphism $D\xra{v}X$ with $D\in Ob(\C{D})$ factors
  {\em uniquely} through $u_X$, i.e.  there exists a {\em unique}
  morphism $D\xra{v'}\app{\C{D}}{X}$ such that $v = u_X\circ v'$.
  Similarly, the co-approximation ${\rm CoApp}_\C{D}(X)$ is the
  approximation of $X$ within $\C{D}^{op}$ viewed as an object of
  $\C{C}^{op}$.  We do not make any assumptions on the existence of
  (co)approximations.
\end{definition}

\begin{theorem}\label{thm:CyclicApproximation}
  Assume $\C{C}$ is a small category with equalizers.  Then the
  approximation $\app{\Lambda}{X_\bullet}$ of any para-(co)cyclic
  module $X_\bullet$ within the category of (co)cyclic modules in
  $\C{C}$ exists.
\end{theorem}

\begin{proof}
  Every para-(co)cyclic object has a canonical endomorphism
  $\omega_\bullet$ defined at each degree $n\geq 0$ by
  $\omega_n\assign \tau_n^{n+1}$ which commutes with all the structure
  morphisms.  The cyclic approximation $\app{\Lambda}{X_\bullet}$ of a
  para-(co)cyclic object $X_\bullet$ is defined degree-wise as the
  equalizer of the pair $(\omega_\bullet,id_\bullet)$ of
  para-(co)cyclic modules in $\C{C}$.  Since both $\omega_\bullet$ and
  $id_\bullet$ are morphisms of para-(co)cyclic module in $\C{C}$,
  their equalizer $\app{\Lambda}{X_\bullet}\xra{} X_\bullet$ is a
  morphism of para-(co)cyclic modules in $\C{C}$.  Moreover,
  $\tau_n^{n+1} = id_n$ on $\app{\Lambda}{X_n}$,
  i.e. $\app{\Lambda}{X_\bullet}$ is a (co)cyclic module in $\C{C}$.
  Assume we have a morphism $Y_\bullet\xra{f_\bullet}X_\bullet$ of
  para-(co)cyclic modules in $\C{C}$ where $Y_\bullet$ is a (co)cyclic
  module in $\C{C}$.  Since $\omega_n f_n = \tau_n^{n+1} f_n = f_n
  \tau_n^{n+1} = f_n$ for any $n\geq 0$, $f_\bullet$ factors through
  the equalizer $\app{\Lambda}{X_\bullet}$.
\end{proof}

\begin{definition}
Recall from \cite{MacLane:Categories} that an endo-functor
$\T{B}\colon\C{C}\xra{}\C{C}$ is called a comonad if there exist
natural transformations $\T{B}\xra{\Delta}\T{B}^2$ and
$\T{B}\xra{\varepsilon}id_\C{C}$ which fit into commutative diagrams
\begin{equation*}
  \xymatrix{
    \T{B}^2(X)  \ar[r]^{\T{B}(\Delta_X)} & \T{B}^3(X)
    & & \T{B}^2(X)  \ar[r]^{\T{B}(\varepsilon_X)} & \T{B}(X) 
        & \ar[l]_{\varepsilon_{\T{B}(X)}} \T{B}^2(X) \\
    \ar[u]^{\Delta_X}\T{B}(X)  \ar[r]_{\Delta_X} & \ar[u]_{\Delta_{\T{B}(X)}} \T{B}^2(X)
    & & & \ar[lu]^{\Delta_X} \T{B}(X) \ar[u]^{id} \ar[ru]_{\Delta_X}
    }
\end{equation*}
However, we will diverge from the standard conventions and we will
refer an object $X$ as a $\T{B}$-comodule if there exists a morphism
$X\xra{\rho_X}\T{B}(X)$ such that the following diagrams commute
\begin{equation*}
  \xymatrix{
    \T{B}(X) \ar[r]^{\Delta_X} & \T{B}^2(X)
    & & \T{B}(X)  \ar[r]^{\varepsilon_X} & X\\
    \ar[u]^{\rho_X} X \ar[r]_{\rho_X} & \T{B}(X) \ar[u]_{\T{B}(\rho_X)}
    & & & \ar[lu]^{\rho_X} X \ar[u]_{id_X}
    }
\end{equation*}
Such objects are called $\T{B}$-coalgebras in
\cite{MacLane:Categories} but later we will work with (co)algebras in
the category of $\T{B}$-comodules and it would have been awkward to
call them ``$\T{B}$-coalgebra coalgebras''.  Also, a morphism
$X\xra{f}Y$ between two $\T{B}$-comodules is called a morphism of
$\T{B}$-comodules if one has a commuting diagram of the form
\begin{equation*}
  \xymatrix{
    \ar[d]_{f} X \ar[r]^{\rho_X} & \T{B}(X) \ar[d]^{\T{B}(f)} \ar[r]^{\varepsilon_X} & X \ar[d]^f \\
    Y \ar[r]_{\rho_Y} & \T{B}(Y) \ar[r]_{\varepsilon_Y} & Y
  }
\end{equation*}
The full subcategory of $\T{B}$-comodules in $\C{C}$ is denoted by
$\C{C}^\T{B}$ and the category of $\T{B}$-comodules and their
morphisms is denoted by $\comod{\T{B}}$.
\end{definition}

\begin{example}
  Let $(\C{C},\otimes)$ be the category of $k$-modules with ordinary
  tensor product of modules as the monoidal product.  Then any
  $k$-coalgebra $(C,\Delta,\varepsilon)$ defines two comonads $(\
  \cdot \ \otimes C)$ and $(C\otimes\ \cdot\ )$.  Moreover, the
  category of comodules in these cases are the same as the category of
  right and left $C$-comodules respectively.
\end{example}

\begin{example}
  Let $(\C{C},\otimes)$ be the opposite category of $k$-modules with
  ordinary tensor product of modules as the monoidal product.  Then
  any $k$-algebra $(A,\mu,e)$ determines two comonads $(\ \cdot\
  \otimes A)$ and $(A\otimes\ \cdot\ )$.  Moreover, the category of
  comodules with respect to these comonads are the same as the
  category of right and left $A$-modules respectively.
\end{example}

\begin{definition}
  A comonad $\T{B}$ is called left exact (resp. right exact) if
  $\T{B}$ commutes with arbitrary small limits (resp. colimits).  In
  other words for any functor $F\colon\C{I}\xra{}\C{C}$ one has
  \[ \lim_\C{I} (\T{B}\circ F) \cong
     \T{B}\left(\lim_\C{I} F\right)
     \text{\ \ \ \ \big(resp.\ } 
     \colim_\C{I} (\T{B}\circ F) \cong
     \T{B}\left(\colim_\C{I} F\right)\big)
  \]
  And a comonad is called exact if it is both left and right exact.
\end{definition}

\begin{definition}
  Let $\T{B}$ be a comonad on a category $\C{C}$.  A para-(co)cyclic
  object $T_\bullet\colon\Lambda_\B{N}\xra{}\C{C}^\T{B}$ is called a
  pseudo-para-(co)cyclic $\T{B}$-comodule if its restriction to the
  subcategory $\Lambda_+$ factors through $\comod{\T{B}}$.
\end{definition}

\begin{theorem}\label{thm:MainApproximation}
  Let $\T{B}$ be a left exact comonad on a complete category $\C{C}$.
  Then every pseudo-para-cyclic $\T{B}$-comodule
  $T_\bullet\colon\Lambda_\B{N}^{op}\xra{}\C{C}^\T{B}$ admits an
  approximation $\app{\Lambda}{T^\T{B}_\bullet}$ within the category
  of cyclic $\T{B}$-comodules.
\end{theorem}

\begin{proof}
  We are going to abuse the notation and use $\partial_j$, $\sigma_i$
  and $\tau_n^\ell$ to denote $T(\partial^n_j)$, $T(\sigma^n_i)$ and
  $T(\tau_n^\ell)$ respectively.  For every $n\geq 0$, denote the
  $\T{B}$-comodule structure morphisms $T_n\xra{}\T{B}(T_n)$ by
  $\rho_n$.  For any $n\geq 0$, define $T^m_n\xra{\eta_{n,m}}T_n$ as
  the equalizer of the pair of morphisms $\T{B}(\tau_n^m) \rho_n$ and
  $\rho_n \tau_n^m$ for every $m\in\B{N}$.  Now define
  \[ T^\T{B}_n\assign \lim_{m\in\B{N}}\ \ T^m_n\xra{\eta_{n,m}}T_n \]
  where $T^\T{B}_n\xra{\eta_n} T_n$ is the canonical morphism into
  $T_n$ for any $n\geq 0$.  Consider the following non-commutative
  diagram in $\C{C}$
  \begin{equation*}
    \begin{CD}
      \T{B}(T_n) @>{\T{B}(\tau_n^j)}>>  \T{B}(T_n)    @>{\T{B}(\tau_n^i)}>>      \T{B}(T_n)\\
      @A{\rho_n}AA                    @AA{\rho_n}A                           @AA{\rho_n}A\\
      T_n        @>>{\ \ \ \tau_n^j\ \ \ }>  T_n   @>>{\ \ \ \tau_n^i\ \ \ }> T_n\\
      @A{\eta_n}AA                      @AA{\eta_n}A\\
      T^\T{B}_n @.                    T^\T{B}_n      
    \end{CD}
  \end{equation*}
  Since $\eta_n$ is the equalizer of the pairs of morphisms
  $\big(\rho_n  \tau_n^i,\ \T{B}(\tau_n^i) \rho_n\big)$
  for all $i\in\B{N}$, if we can show that
  \begin{equation}\label{eq:Approximation1}
    \rho_n  \tau_n^i  \tau_n^j \eta_n
    = \T{B}(\tau_n^i) \rho_n  \tau_n^j \eta_n
  \end{equation}
  for all $i\in\B{N}$ we will obtain a functorial `restriction' of
  $\tau_n^j$ to $T^\T{B}_n$ which will be denoted by
  $(\tau_n^j)^\T{B}$ for any $j\in\B{N}$.  Consider the left hand side
  of Equation~\ref{eq:Approximation1} which is
  \begin{align*}
    \rho_n  \tau_n^{i+j} \eta_n
    = \T{B}(\tau_n^{i+j}) \rho_n \eta_n
    = \T{B}(\tau_n^i)  \T{B}(\tau_n^j) \rho_n \eta_n
    = \T{B}(\tau_n^i)  \rho_n  \tau_n^j \eta_n
  \end{align*}
  as we wanted to show.

  Now, for $0\leq j\leq n+1$ consider the following diagram in $\C{C}$
  \begin{equation*}
    \begin{CD}
      \T{B}(T_{n+1}) @>{\T{B}(\partial_j)}>>  \T{B}(T_n)    @>{\T{B}(\tau_n^i)}>>  \T{B}(T_n)\\
      @A{\rho_{n+1}}AA                        @AA{\rho_n}A                        @AA{\rho_n}A\\
      T_{n+1}     @>>{\ \ \ \partial_j\ \ \ }>  T_n      @>>{\ \ \ \tau_n^i\ \ \ }> T_n\\
      @A{\eta_{n+1}}AA                          @AA{\eta_n}A\\
      T^\T{B}_{n+1} @.                           T^\T{B}_n      
    \end{CD}
  \end{equation*}
  where the square on top right does not commute and square on top
  left commutes as long as $0\leq j\leq n$.  However, since
  $\partial_{n+1} = \partial_0 \tau_{n+1}$ (recall that $T_\bullet$ is
  cyclic not cocyclic) and $\tau_{n+1}$ has a restriction to
  $T^\T{B}_{n+1}$, one can assume WLOG that $0\leq j\leq n$.  If we
  can show that
  \begin{equation}\label{eq:Approximation}
    \rho_n  \tau_n^i  \partial_j \eta_{n+1}  
     = \T{B}(\tau_n^i) \rho_n  \partial_j \eta_{n+1}  
  \end{equation}
  for any $i\in\B{N}$, one obtains a unique morphism
  $T^\T{B}_{n+1}\xra{}T^\T{B}_n$ which is going to be denoted as
  $(\partial_j)^\T{B}$.  The uniqueness of this morphism implies its
  functoriality.  Consider the left hand side of the
  Equation~\ref{eq:Approximation}
  \begin{align*}
    \rho_n  \tau_n^i  \partial_j \eta_{n+1}  
     = & \rho_n  \partial_q  \tau_{n+1}^{i+p}  \eta_{n+1}
  \end{align*}
  where $(i+j)=(n+1)p+q$ and $0\leq q\leq n$.  Now use the fact that
  $0\leq j\leq n$ and $T_\bullet$ is a pseudo-para-cyclic to deduce
  \begin{align*}
    \rho_n  \partial_q  \tau_{n+1}^{i+p}  \eta_{n+1}
    = & \T{B}(\partial_q) \rho_n  \tau_{n+1}^{i+p}  \eta_{n+1}
    = \T{B}(\partial_q)  \T{B}(\tau_{n+1}^{i+p}) \rho_n  \eta_{n+1}\\
    = & \T{B}(\tau_{n+1}^i) \T{B}(\partial_j) \rho_n  \eta_{n+1}
    =  \T{B}(\tau_{n+1}^i) \rho_n  \partial_j \eta_{n+1}
  \end{align*}
  as we wanted to show.  One can similarly prove that the relevant
  diagrams commute for the degeneracy morphisms.  This finishes the
  proof that $T^\T{B}_\bullet$ is a para-cyclic module in $\C{C}$.

  Now, for an arbitrary $j\in\B{N}$ consider the non-commutative
  diagram
  \begin{equation*}
    \begin{CD}
      T_n    @>{\rho_n}>>  \T{B}(T_n)  @>{\T{B}(\rho_n)}>>   \T{B}^2(T_n)\\
      @A{\tau_n^j}AA         @AA{\T{B}(\tau_n^j)}A               @AA{\T{B}^2(\tau_n^j)}A\\
      T_n   @>>{\rho_n}>   \T{B}(T_n)  @>>{\T{B}(\rho_n)}>   \T{B}^2(T_n)\\
      @A{\eta_n}AA           @AA{\T{B}(\eta_n)}A\\
      T^\T{B}_n    @.        \T{B}(T^\T{B}_n)
    \end{CD}
  \end{equation*}
  and the composition
  \begin{align*}
    \T{B}^2(\tau_n^j)  \T{B}(\rho_n) \rho_n \eta_n
     = & \T{B}^2(\tau_n^j)  \Delta_{T_n} \rho_n \eta_n
     =  \Delta_{T_n}  \T{B}(\tau_n^j) \rho_n \eta_n\\
     = & \Delta_{T_n} \rho_n  \tau_n^j \eta_n
     =  \T{B}(\rho_n) \rho_n  \tau_n^j \eta_n\\
     = & \T{B}(\rho_n)  \T{B}(\tau_n^j) \rho_n \eta_n
  \end{align*}
  The equality of the first and the last terms implies $\rho_n \eta_n$
  factors through the limit of the equalizers of the pairs
  $\T{B}^2(\tau_n^j) \T{B}(\rho_n)$ and $\T{B}(\rho_n)
  \T{B}(\tau_n^j)$ as $j$ runs through the set of all natural numbers.
  But $\T{B}$ is a left exact comonad which means this limit is
  exactly $\T{B}(T^\T{B}_n)$.  Thus we get the $\T{B}$-comodule
  structure on $T^\T{B}_n$ which implies $T^\T{B}_\bullet$ is a
  para-cyclic module in $\C{C}^\T{B}$.

  Now we need to show given any morphism $[n]\xra{\phi}[m]$ in
  $\Lambda_\B{N}$ its image $\phi^\T{B}$ under the newly constructed
  functor $T^\T{B}_\bullet$ is a morphism of $\T{B}$-comodules.  In
  order to prove this fact we need the following diagram to commute
  \begin{equation*}
    \begin{CD}
      \T{B}(T^\T{B}_n)  @<{\T{B}(\phi)}<< \T{B}(T^\T{B}_m)\\
      @A{\rho_n}AA                      @AA{\rho_m}A\\
      T^\T{B}_n         @<<{\phi^\T{B}}<  T^\T{B}_m
    \end{CD}
  \end{equation*}
  To achieve this, first we need to show that the larger squares in
  the following diagrams commute
  \begin{equation*}
    \begin{CD}
      \T{B}(T_n) @<{\T{B}(\partial_j)} << \T{B}(T_{n+1})
      @.\hspace{2cm} @. \T{B}(T_n) @>{\T{B}(\sigma_j)}>> \T{B}(T_{n+1})\\
      @A{\rho_n}AA      @AA{\rho_n}A        @.   @A{\rho_n}AA   @AA{\rho_m}A\\
      T_n   @<{\partial_j}<< T_{n+1}
      @.\hspace{2cm} @. T_n   @>{\sigma_j}>>   T_{n+1}\\
      @A{\eta_n}AA     @AA{\eta_n}A @.   @A{\eta_n}AA        @AA{\eta_n}A\\
      T^\T{B}_n  @<<{(\partial_j)^\T{B}}<  T^\T{B}_{n+1}
      @.\hspace{2cm} @. T^\T{B}_n  @>>{(\sigma_j)^\T{B}}>  T^\T{B}_{n+1}
    \end{CD}
  \end{equation*}
  for any $i\geq 0$ and $0\leq j\leq n$.  In these diagrams, the top
  squares commute since $T_\bullet$ is pseudo-para-cyclic.  We already
  have shown the bottom squares commute.  Thus both diagrams commute
  for the prescribed range.  Then we must show that the larger square
  in the following diagram commutes
  \begin{equation*}
    \begin{CD}
    \T{B}(T_n) @>{\T{B}(\tau_n^i)}>> \T{B}(T_n) \\
    @A{\rho_n}AA      @AA{\rho_n}A\\
    T_n   @>>{\tau_n^i}> T_n  \\
    @A{\eta_n}AA      @AA{\eta_n}A\\
    T^\T{B}_n  @>>{(\tau_n^i)^\T{B}}> T^\T{B}_n
    \end{CD}
  \end{equation*}
  The bottom square commutes while the top square does not.  However,
  $\eta_n$ equalizes $\rho_n  \tau_n^i$ and
  $\T{B}(\tau_n^i) \rho_n$.  Therefore the larger diagram
  commutes.  This finally finishes the proof that $T^\T{B}_\bullet$ is
  a para-cyclic $\T{B}$-comodule.  Now we use
  Theorem~\ref{thm:CyclicApproximation} to finish the proof.
\end{proof}

\begin{remark}
  There are 8 versions of Theorem~\ref{thm:MainApproximation}
  \begin{quote}
    A pseudo-para-(co)cyclic $\T{B}$-(co)module in $\C{C}$ admits
    an(a) (co)approximation in the category of (co)cyclic
    $\T{B}$-(co)modules.
  \end{quote}
  However, since the proof is given for an arbitrary complete
  category, by assuming $\C{C}$ is both complete and cocomplete, one
  can use $\C{C}$ and $\C{C}^{op}$ interchangibly.  This reduces the
  number of versions to 4:
  \begin{quote}
    A pseudo-para-(co)cyclic $\T{B}$-comodule in $\C{C}$ admits an(a)
    (co)approximation in the category of (co)cyclic $\T{B}$-comodules.
  \end{quote}
  From the remaining 3, we are interested in the following
\end{remark}

\begin{theorem}
  Let $\T{B}$ be a left exact comonad on a complete category $\C{C}$.
  Then every pseudo-para-cocyclic $\T{B}$-comodule
  $T_\bullet\colon\Lambda_\B{N}\xra{}\C{C}^\T{B}$ admits an
  approximation $\app{\Lambda}{T^\T{B}_\bullet}$ within the category
  of cocyclic $\T{B}$-comodules.
\end{theorem}

\begin{proof}
  As before let $\rho_n$ denote the $\T{B}$-comodule structure
  morphism on $T_n$ for any $n\geq 0$.  Let $\Gamma(n)$ be the set of
  pairs of morphism of the form
  \[  \big(\rho_n  \tau_n^i,\ \T{B}(\tau_n^i) \rho_n\big)
      \text{\ for $i\geq 0$\ \ \ \ or\ \ \ \ }
      \big(\rho_n  \partial_{n+1},\ \T{B}(\partial_{n+1}) \rho_n\big)
   \]    
  and we define $T^\gamma_n\xra{\eta(\gamma)}T_n$ as the equalizer of
  a pair $\gamma\in\Gamma(n)$.  Next we define the approximation
  $T^\T{B}_n$ for each $n\geq 0$ as
  \begin{equation*}
    T^\T{B}_n\assign \lim_{\gamma}\ T^\gamma_n\xra{\eta(\gamma)}T_n
  \end{equation*}
  where we use $T^\T{B}_n\xra{\eta_n}T_n$ to denote the canonical
  morphism into $T_n$.  The rest of the proof is very similar to that
  of Theorem~\ref{thm:MainApproximation} and we leave it to the
  interested reader to finish it.
\end{proof}

\begin{definition}
  The (co)cyclic $\T{B}$-comodule $\app{\Lambda}{T_\bullet^\T{B}}$
  corresponding to a pseudo-para-(co)cyclic $\T{B}$-comodule
  $T_\bullet$ is called the universal (co)cyclic $\T{B}$-comodule of
  $T_\bullet$.  Moreover, given a functor of the form $\C{F}\colon
  \C{C}^\T{B}\xra{}\lmod{k}$ and a (co)homology functor $\C{H}_*$ on
  the category of (co)cyclic $k$-modules, one can compute
  \[  \C{H}_* \C{F}(\app{\Lambda}{T_\bullet^\T{B}}) \]
  We will call this (co)homology as {\em the $\T{B}$-equivariant
  $\C{H}$-(co)homology of $T_\bullet$ with coefficients in $\C{F}$.}
\end{definition}

\section{The universal cyclic theory of (co)module (co)algebras}
\label{sec:Examples}

\subsection{Hopf and equivariant cyclic theory of module coalgebras}

Fix a commutative unital ring $k$ and an associative/ coassociative
unital/counital $k$-bialgebra $(B,\mu_B,\B{I},\Delta_B,\varepsilon)$.
Our base category is the opposite category of $k$-modules with the
opposite tensor product over $k$,
i.e. $(\C{C},\otimes)\assign(\lmod{k}^{op}, \otimes_k^{op})$.  Our
base comonad in $\C{C}$ is going to be $\T{B}\assign (\ \cdot \
\otimes B)$.  Since we defined the comonad in the opposite category,
we will use the algebra structure on $B$.

The category of left $B$-modules (i.e. $\T{B}$-comodules in $\C{C}$)
is a monoidal category with respect to the ordinary tensor product of
$k$-modules with the diagonal action of $B$ on the left.  Explicitly,
given a pair of $B$-modules $X$ and $Y$, the $B$-module structure on
the product is given by
\[ b(x\otimes y) \assign b_{(1)}x\otimes b_{(2)}y
\]
for any $x\otimes y\in X\otimes Y$.  However, the product is not
symmetric unless $B$ is cocommutative but there is a braided monoidal
structure if one restricts oneself to use Yetter-Drinfeld modules.  If
denote the full subcategory of left $B$-modules of $\lmod{k}$ by
$\C{L}(B)$ then one can see that $\C{C}^\T{B} = \C{L}(B)^{op}$.

Fix a left/left $B$-module/comodule $M$ and for each $X\in
Ob(\C{C}^\T{B})$ define a transposition $w_{M,X}:M\otimes
X\xra{}X\otimes M$ by
\[ w_{M,X}(m\otimes x)\assign m_{(-1)}x\otimes m_{(0)} \]
for any $m\otimes x\in M\otimes X$.  

Any algebra $(C,\delta_C,\varepsilon)$ in $\C{C}^\T{B}$ is a
$B$-module coalgebra and therefore is automatically $w$-transpositive.
We form the objects $P_\bullet(C,M)$ and $T_\bullet(C,M)\assign
\colim_S P_\bullet(C,M)$ in $(\lmod{k}^{op},\otimes_k^{op})$ and
consider the latter as a para-cyclic module in $\C{C}$.  In fact
$T_\bullet(C,M)$ carries a pseudo-para-cyclic $\T{B}$-module structure
and
\[ Q_\bullet(C,M)^{op} = {\rm App}_\Lambda(T_\bullet(C,M)^\T{B}) \]
where $Q_\bullet(C,M)$ is the Hopf-equivariant cocyclic object defined
in \cite{Kaygun:BivariantHopf} for a $B$-module coalgebra $C$ and an
arbitrary $B$-module/comodule $M$.  Therefore, the Hopf cyclic
(co)homology of the triple $(C,B,M)$ is defined as the cyclic
(co)homology of the cocyclic $k$-module $C_\bullet(C,M)\assign
k\otimes_{B}Q_\bullet(C,M)$.  

Similarly, any coalgebra $(A,\mu_A,1)$ in $\C{C}^\T{B}$ is a
$B$-module algebra and therefore is automatically $w$-transpositive.
We form the objects $P_\bullet(A,M)$ and $T_\bullet(A,M)\assign
\colim_S P_\bullet(A,M)$ in $(\lmod{k}^{op},\otimes_k^{op})$ and
consider the latter as a para-cocyclic module in $\C{C}$.  In fact
$T_\bullet(A,M)$ carries a pseudo-para-cocyclic $\T{B}$-module
structure and we see that
\[ Q_\bullet(A,M)^{op} = {\rm App}_\Lambda(T_\bullet(A,M)^\T{B}) \]
where $Q_\bullet(A,M)$ is the Hopf-equivariant cyclic object defined
in \cite{Kaygun:BivariantHopf} for a $B$-module algebra $A$ and an
arbitrary $B$-module/comodule $M$.  Therefore, the Hopf cyclic
(co)homology of the triple $(A,B,M)$ is defined as the cyclic
(co)homology of the cyclic $k$-module $C_\bullet(A,M)\assign
k\otimes_{B}Q_\bullet(A,M)$.

\subsection{Hopf-Hochschild homology}

Let $(\C{C},\otimes)$, $M$, $w$, $B$ and $\T{B}$ be as before.  Assume
$A$ is a $B$-module algebra and construct $P_\bullet(A,M)$ in $\C{C}$.
However, since $A$ is a $B$-module algebra, $P_\bullet(A,M)$ is a
functor from $S$ into $\mathcal{C}^\T{B}$, i.e. it is a graded
$B$-module.  This time, instead of considering colimit of
$P_\bullet(A,M)\to \mathcal{C}$, we consider the colimit of
$P_\bullet(A,M)\colon S\to \C{C}^\T{B}$, which we denote by
$T'_\bullet(A,M)$.  This colimit is a little more than
pseudo-para-cyclic $B$-module: viewed just as a simplicial object, it
is actually a simplicial $B$-module.  The Hochschild homology of
$k\otimes_B T'_\bullet(A,M)$ is the Hopf--Hochschild homology of the
triple $(A,B,M)$ as constructed in \cite{Kaygun:HopfHochschild}.

\subsection{Hopf and equivariant cyclic theory of comodule (co)algebras}

Fix a commutative unital ring $k$ and an associative/ coassociative
unital/counital $k$-bialgebra $(B,\mu_B,\B{I},\Delta_B,\varepsilon)$.
Our base category is the category of $k$-modules with the ordinary
tensor product over $k$, i.e. $(\C{C},\otimes)\assign(\lmod{k},
\otimes_k)$.  Our base comonad in $\C{C}$ is going to be $\T{B}\assign
(B\otimes \ \cdot \ )$ thus we will use the coalgebra structure on
$B$.

The category of left $B$-comodules (i.e. $\T{B}$-comodules in $\C{C}$)
is a monoidal category with respect to the ordinary tensor product of
$k$-modules with the diagonal coaction of $B$ on the left.
Explicitly, given a pair of $B$-modules $X$ and $Y$, the $B$-comodule
structure on the product is given by
\[ \rho(x\otimes y) \assign x_{(-1)}y_{(-1)}\otimes (x_{(0)}\otimes y_{(0)})
\]
for any $x\otimes y\in X\otimes Y$.  However, the product is not
symmetric unless $B$ is commutative but there is a braided monoidal
structure if one restricts oneself to use Yetter-Drinfeld modules.

Fix a left/left $B$-module/comodule $M$ and for each $X\in
Ob(\C{C}^\T{B})$ define a transposition $w_{M,X}:M\otimes
X\xra{}X\otimes M$ by
\[ w_{M,X}(m\otimes x)\assign x_{(0)}\otimes x_{(-1)}m \]
for any $m\otimes x\in M\otimes X$.  

Any coalgebra $(C,\delta_C,\varepsilon)$ in $\C{C}^\T{B}$ is a
$B$-comodule coalgebra and therefore is automatically $w$-transpositive.
We form the objects $P_\bullet(C,M)$ and $T_\bullet(C,M)\assign
\colim_S P_\bullet(C,M)$ in $(\lmod{k},\otimes_k)$ and consider the
latter as a para-cocyclic module in $\C{C}$.  In fact $T_\bullet(C,M)$
carries a pseudo-para-cocyclic $\T{B}$-comodule structure and
\[ Q_\bullet(C,M) \assign {\rm App}_\Lambda(T_\bullet(C,M)^\T{B}) \]
is a cocyclic $\T{B}$-comodule, i.e. a cocyclic $B$-comodule.

Similarly, any algebra $(A,\mu_A,1)$ in the category $\C{C}^\T{B}$ is
a $B$-comodule algebra and therefore is automatically $w$-transpositive.
We form the objects $P_\bullet(A,M)$ and $T_\bullet(A,M)\assign
\colim_S P_\bullet(A,M)$ in $(\lmod{k},\otimes_k)$ and consider the
latter as a para-cyclic module in $\C{C}$.  In fact $T_\bullet(A,M)$
carries a pseudo-para-cyclic $\T{B}$-comodule structure and we see that
\[ Q_\bullet(A,M) \assign {\rm App}_\Lambda(T_\bullet(A,M)^\T{B}) \]
is a cyclic $\T{B}$-comodule, i.e. a cyclic $B$-comodule.  Moreover,
the cyclic cohomology of the cyclic $k$-module $k\otimes_B
Q_\bullet(A,M)$ is the bialgebra cyclic homology of a module coalgebra
as defined in ~\cite{Kaygun:BialgebraCyclicK}.


\end{document}